\documentclass{article}
\usepackage{latexsym}
\usepackage{bm}
\usepackage{amssymb}
\usepackage{amsthm}
\usepackage{amsmath}
\usepackage{graphicx}
\usepackage{calrsfs}
\usepackage[all]{xy}
\usepackage[linesnumbered,boxed]{algorithm2e}
\usepackage{float}
\usepackage[usenames,dvipsnames]{color}
\usepackage{listings}
\usepackage{textcomp}
\usepackage[T1]{fontenc}
\usepackage[Q=yes]{examplep}
\usepackage{courier}
\usepackage{subfigure}% in preamble

\definecolor{listinggray}{gray}{0.9}
\definecolor{lbcolor}{rgb}{0.9,0.9,0.9}
\begin{document}
  % set up commands
  \newcommand{\rsup}[1]{^{(#1)}} % round bracket superscript
  \newcommand{\Poincare}{Poincar\'e\xspace}
  \newcommand*{\rom}[1]{\expandafter\@slowromancap\romannumeral #1@}
  \newcommand{\setr}{\mathbb{R}}
  \newcommand{\supp}{\mathrm{supp}}
%  \newcommand{\span}{\mathrm{span}}
  % set up environment
  \newtheorem{theorem}{Theorem}[section]
  \newtheorem{lemma}[theorem]{Lemma}
  \newtheorem{proposition}[theorem]{Proposition}
  \newtheorem{corollary}[theorem]{Corollary}
  \newtheorem{assumption}[theorem]{Assumption}
  \newenvironment{myproof}[1][Proof.]{\begin{trivlist}
  \item[\hskip \labelsep {\bfseries #1}]}{\end{trivlist}}
  \newenvironment{definition}[1][Definition]{\begin{trivlist}
  \item[\hskip \labelsep {\bfseries #1}]}{\end{trivlist}}
  \newenvironment{example}[1][Example]{\begin{trivlist}
  \item[\hskip \labelsep {\bfseries #1}]}{\end{trivlist}}
  \newenvironment{remark}[1][Remark]{\begin{trivlist}
  \item[\hskip \labelsep {\bfseries #1}]}{\end{trivlist}}
  % title
  \title{Sparse Recovery with Oversampling Ratio Greater than One Half}
  \author{Wenbin Zhang}
  \maketitle
%  \tableofcontents
  
%\section{Notations}
%\begin{itemize}
%\item $x\in \setr^n$ signal, $A\in \setr^{m\times n}$ sensing matrix, $b\in \setr^m$ measurements. ($m\leq n$)
%\item $supp(x) = S$, $|S|\leq s$.
%\item $i\in \{1,\ldots, N\}$: row index; $j\in \{1,\ldots, N\}$: column index.
%\item Superscript of $x$ means the type of $x$: plain $x$ means true solution; $x'$, $x''$ denote other solutions if they exist; $k$-th iterate of $x$ is denoted by $x^{(k)}$.
%\item Subscript of $x$ means some entry or some entries of $x$: $x_j$ denotes $j$-th entry of $x$. $x_T, T\in N$ has two meanings: one is the truncation of $x$ on support $T$, the other is zeroing out $x$ on indexes other than $T$. 
%\item $S$, $T$, $S\rsup{n}, T\rsup{n}\in \{1,\ldots, N\}$: some support names.
%\item $\Phi, \Psi, \Phi(S), \Phi(S', S'')\subseteq \setr^N$: some sets of interest related to $x$.
%\item $v, z\in \setr^N$: some vectors of interest related to $x$.
%\item $s,t,r\in \mathbb{N}$: support length.
%\item $k,l,n,p,q\in \mathbb{N}$: iteration indexes.
%\end{itemize}

\section{Introduction}
Sparse recovery technique (or compressive sensing ~\cite{candes2006robust,donoho2006compressed} technique) offers to recover signals under sparsity prior with fewer measurements than Nyquist compression would require. Mathematically, the standard sparse recovery problem can be formulated as the following:
\begin{equation}
(P_{x})\quad
  \begin{cases}
  \quad \min_{z\in \setr^N} & f(z) = \frac{1}{2}\|Az-b\|_2^2 \nonumber\\
  \quad \mbox{s.t. } & \|z\|_0\leq s, \nonumber
  \end{cases}
\end{equation}
where sensing matrix $A\in \setr^{m\times N}$, sparse signal $x\in \setr^N(\|x\|_0\leq s)$, dummy variable $z\in \setr^N$, and measurement vector $b = Ax \in \setr^m$. Various algorithms have been developed to tackle $(P_x)$ or varieties of $(P_x)$ such as orthogonal matching pursuit~\cite{pati1993orthogonal}, iterative hard thresholding~\cite{blumensath2009iterative}, CoSaMP~\cite{needell2009cosamp}, hard thresholding pursuit~\cite{foucart2011hard}, subspace pursuit~\cite{dai2009subspace}, etc. Successful sparse recovery algorithms possess a list of properties~\cite{needell2009cosamp} including minimal number of measurements required, adaptability to different sampling schemes, noise robusty, optimal error guarantees, and efficient computational resource usage. 

In this paper, however, we are interested in designing algorithms that can successfully solve $(P_x)$ with number $m$ of measurements as small as possible. This leads to the solvability of $\ell_0$ sparse recovery problem $(P_x)$ with respect to oversampling ratio $\rho=s/m$~\cite{donoho2004neighborly,donoho2006high}. Unfortunately, if we require $(P_x)$ to be solvable for ALL $\|x\|_0\leq s$, then there is a theoretic upper-bound for oversampling ratio $\rho$. which is stated in the following theorem:

\begin{theorem}\label{thm_m>=2s}
Given a sensing matrix $A\in \setr^{m\times N}$ ($m\leq N$). The sparse recovery problem ($P_{x}$) has unique solution for all $\|x\|_0\leq s$ if and only if all the $m$-by-$2s$ sub-matrices of $A$ have full column ranks, implying that the oversampling ratio $s/m\leq 1/2$.
\end{theorem}
\begin{myproof}
See \cite{foucart2013mathematical}, Theorem 2.13.\qed
\end{myproof}

\begin{figure}[!h]
\centering 
\subfigure[HTP phase transition]{%
\includegraphics[scale=0.30]{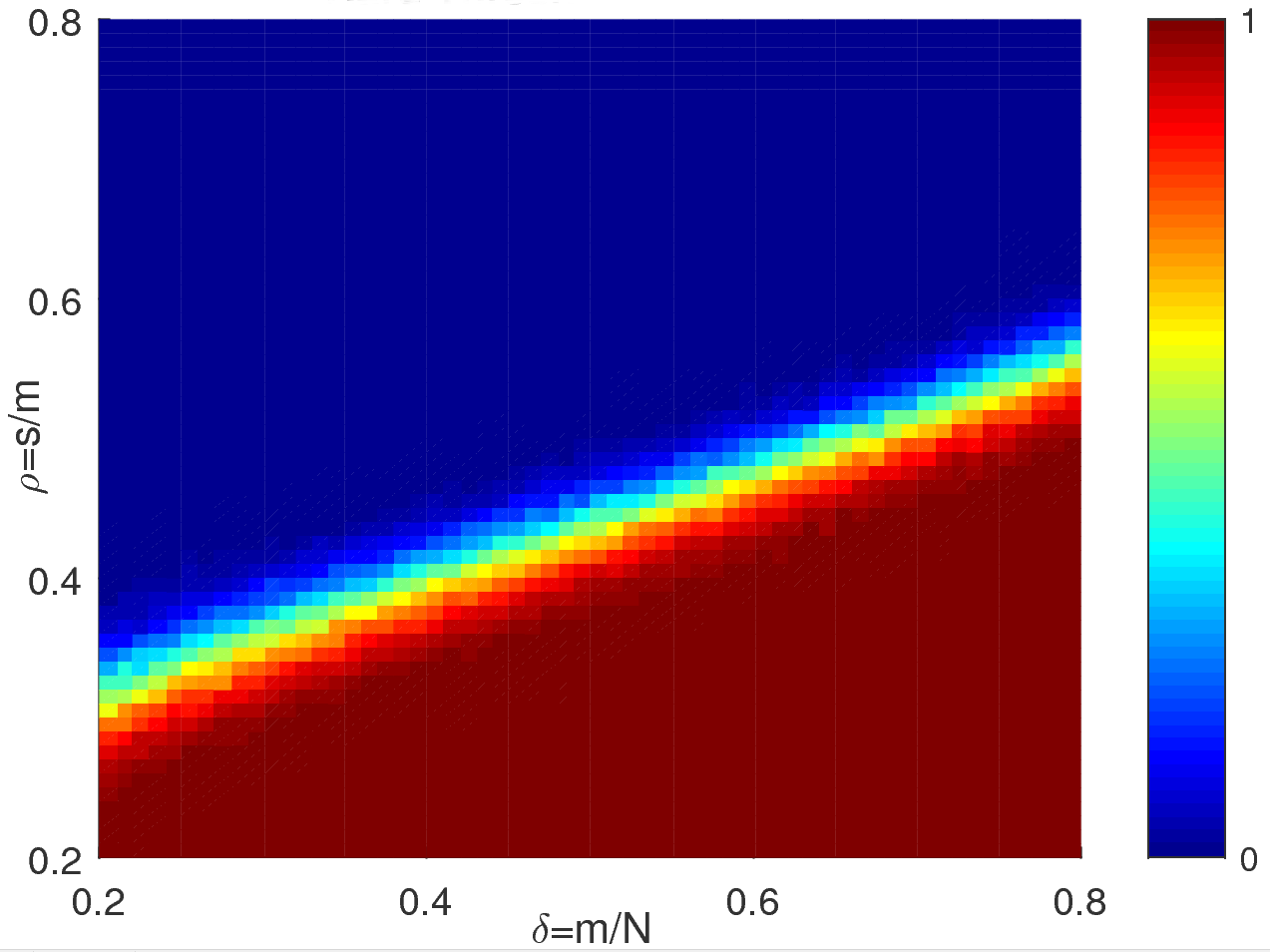}
\label{fig:phase_HTP}}
\quad
\subfigure[DP-HTP phase transition]{%
\includegraphics[scale=0.30]{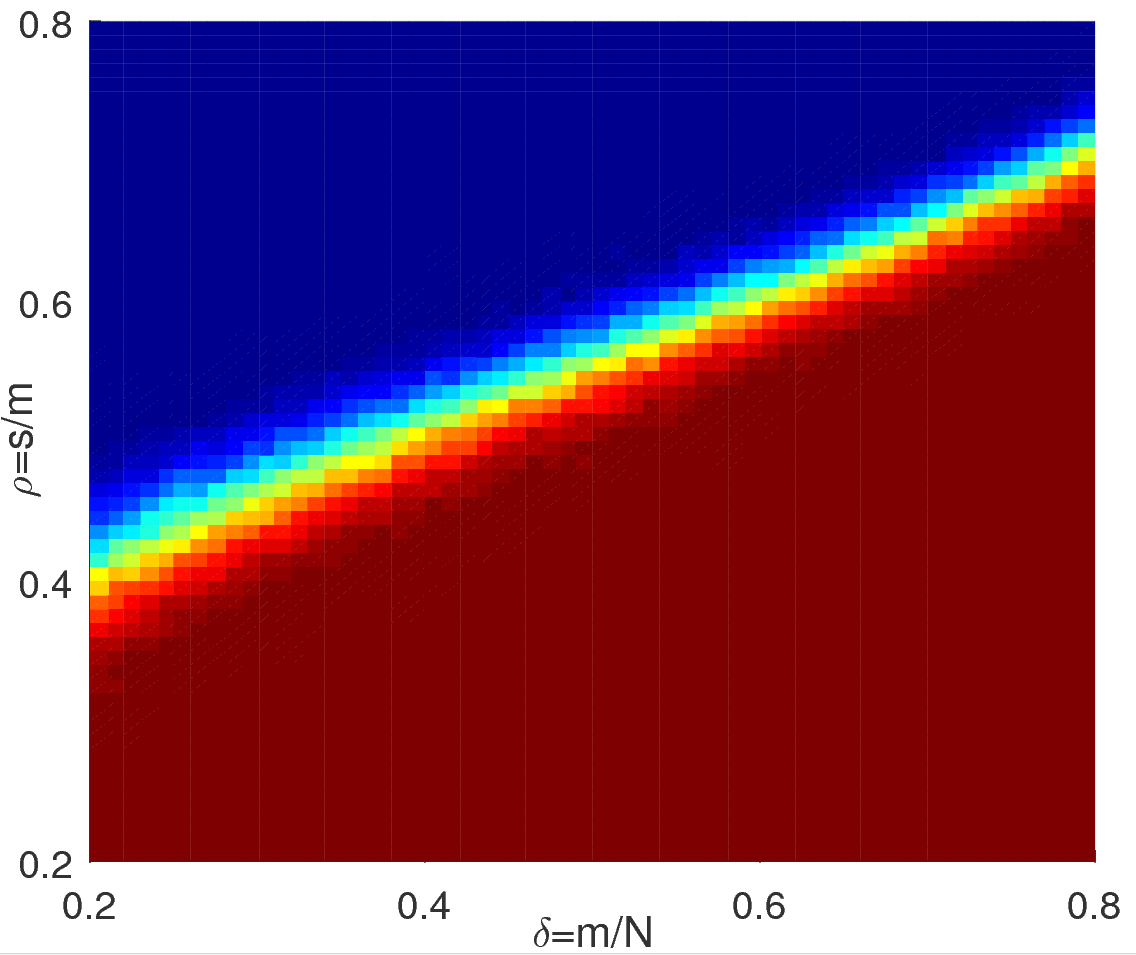}
\label{fig:phase_DPHTP}}

\caption{phase transition: HTP vs DP-HTP}
\label{fig:phase_HTPvsDPHTP}
\end{figure}

However, we observe that the phase transition behaviour of a certain modified Hard Thresholding Pursuit algorithm (named DP-HTP) seems to contradict theorem \ref{thm_m>=2s}, see Fig. \ref{fig:phase_HTPvsDPHTP}. Note that the large area in the phase transition map of DP-HTP with oversampling ratio $s/m\geq 1/2$ means DP-HTP algorithm successfully recovers the true sparse signals.

Here, we attempt to answer the following two questions:
\begin{enumerate}
\item[1] Is it possible for $\ell_0$ sparse recovery problem $(P_x)$ to have unique solution, in some sense, with oversampling ratio $s/m$ larger than one half?
\item[2] How to explain the phase transition behaviour of DP-HTP algorithm?
\end{enumerate}

The remainder of the paper is organized as follows. Section 2 gives a somehow positive answer to question 1 as well as a full theoretical result on uniqueness of solution to $(P_x)$ with respect to oversampling ratio $\rho=s/m$. Section 3 presents the algorithm DP-HTP and justifies its phase-transition behaviour. Section 4 presents extended numerical experiments of HTP and DP-HTP with respect to different sparse signal $x$'s as well as different sensing matrix $A$s.
%\subsection{Map of the article}

\section{Uniqueness of solution when oversampling ratio greater than one half}
The existence of solution to $(P_x)$ is guaranteed since we generate this problem using a sparse vector $x$ and thus $x$ itself must be a solution. Our concern is the uniqueness.
 
It is trivial to observe that when oversampling ratio $\rho=s/m>1$, the solution to $(P_x)$ is never unique since sub-matrix $A_S\in \setr^{m\times s}$ has more columns than rows, where index set $S$ denotes the support of the true sparse signal $x$, and even if we know the true support, there are still infinitely many solutions to $(P_x)$.

To prevent the above situation from happening, we introduce a basic assumption on the sensing matrix $A\in \setr^{m\times N}$:

\begin{assumption}\label{basic_assum}
$A_T$ is non-singular for any index set $T\subset \{1,\ldots, N\}$ and $|T|\leq m$.
\end{assumption}

Conditions on sub-matrix of the sensing matrix is quite common in compressive sensing community, such as the restricted isometry property(RIP)~\cite{candes2005decoding}.\textbf{Assumption \ref{basic_assum}} is just a weaker version of RIP. Also note that sensing matrix $A$ has overwhelming large probability to satisfy the basic \textbf{Assumption \ref{basic_assum}} if $A$ is generated using usual random generating techniques such as random gaussian or random DCT ~\cite{foucart2013mathematical}.

Since we already have negative answer for the uniqueness when $\rho=s/m>1$, the rest situations to be discussed are situation $1/2<\rho<1$ and $\rho=1$.

\subsection{Situation $1/2<\rho<1$}
We characterize the following set, which is the set of $s$-sparse $x$ such that ($P_x$) has multiple solutions:
\begin{equation}
\Phi = \left\{x\in \setr^N \big |\  \|x\|_0 \leq s, \mbox{($P_x$) has multiple solutions}\right\}.
\end{equation}
\begin{theorem}\label{thm_m<2s}
Given a fixed sensing matrix $A\in \setr^{m\times N}$ and we suppose the basic \textbf{Assumption \ref{basic_assum}} holds on $A$. If $1/2<\rho<1$, or $s<m<2s$, then the non-uniqueness set $\Phi$ has the following decomposition:
\begin{equation}
\Phi = \bigcup_{\begin{subarray}{l}
                S',S''\subseteq \{1,\ldots, N\}\\|S'|,|S''|\leq s
                \end{subarray}} 
       \Phi (S', S''),
\end{equation}
where 
\begin{align}
&\Phi(S',S'') =  \nonumber\\
               &\qquad \big\{x'\in \setr^N \big | \supp(x)=S',\exists x''\neq x', \supp(x'')=S'', A(x'-x'')=0 \big\},
\end{align}
and 
\begin{equation}
\dim\big(\mathrm{span}(\Phi(S',S''))\big)\leq 2s-m <s.
\end{equation}
\end{theorem}

\begin{myproof}
We decompose $\Phi$ first to the support of $x'$ and then further to the support of another solution:
\begin{align}
\Phi =& \bigcup_{\begin{subarray}{l}
                S'\subseteq N, |S'|\leq s
                \end{subarray}} 
       \big\{
         x'\in \setr^N \big |\  \supp(x') = S', \mbox{($P_{x'}$) has multiple solutions.}
       \big\}\\
     =& \bigcup_{\begin{subarray}{l}
                S',S''\subseteq N\\|S'|,|S''|\leq s
                \end{subarray}} 
        \Phi (S', S'').
\end{align}
For any $x'\in \Phi(S',S'')$, there exists an $x''\neq x'$ with $\supp(x'')=S''$ and 
\begin{equation*}
A(x'-x'')\equiv Av=0.
\end{equation*}
We first let $s'=|S'|$, $s''=|S''|$, $T=S'\cup S''$, $l=|S'\cap S''|$ and then analyze the constraints imposed on $x'$. We will only discuss indexes in support $T$ and this may help to eliminate the ambiguity of some symbols.

Of course, $|T|=s'+s''-l$ and the existence of $x''$ implies 
\begin{equation}
\dim(\ker(A_T))=s'+s''-l-m >0.
\end{equation}
The vector $v_T$ has a special structure in regard to $x'$ and $x''$:
\begin{equation}
(v_T)_j=
\begin{cases}
x'_j         &\mbox{if } j\in S'\backslash S'' \\
x'_j- x''_j  &\mbox{if } j\in S'\cup S'' \\
-x''_j       &\mbox{if } j\in S''\backslash S',
\end{cases}
\end{equation}
and accordingly, we separate $x'$ into two parts: $x'_{S'\backslash S''}$ and $x'_{S'\cup S''}$.

For index set $S'\backslash S''$, we let 
\begin{align}
&\Psi(S',S'') = \nonumber\\ 
   &\quad\quad\big\{ 
                z_{S'\backslash S''}\in \setr^N \big |\ z\in \setr^N, \supp(z)=T, z_T \in \ker(A_T), T=S'\cup                             
                S''
               \big\}.
\end{align}
Since $\Psi(S',S'')$ is some projection image of $\ker(A_T)$, its dimension is no larger than that of $\ker(A_T)$:
\begin{equation}
\dim(\Psi(S',S'')) \leq \dim(\ker(A_T)).
\end{equation}
Of course, $x'_{S'\backslash S''}\in \Psi(S',S'')$.

Now we consider index set $S'\cap S''$. Although we have constraints on $x'-x''$, $x'$ itself on the joint support $S'\cap S''$ is actually arbitrary. Thus we have
\begin{equation}
x'_{S'\cap S''}\in \Omega(S',S'') = \big\{z_{S'\cap S''}\big |\ z\in \setr^N\big\},
\end{equation}
of which the dimension is just $l$: $\dim(\Omega(S',S''))=l$.

Since the full $x'=x'_{S'\backslash S''} + x'_{S'\cup S''}$, we have 
\begin{equation}
x'\in \Psi(S',S'') \oplus \Omega(S',S''),
\end{equation}
and
\begin{align}
\dim(\Phi(S',S''))& \leq \dim(\Psi(S',S'')\oplus \Omega(S',S'')) \nonumber\\
                  & \leq s'+s''-l-m+l=s'+s''-m\leq 2s-m.
\end{align}
As long as $s<m$, the dimension of $\Phi(S',S'')$ is less than $s$, as is stated in the theorem.
\qed
\end{myproof}

\begin{remark}
Those ``bad'' $x$'s are contained in a finite union of subspaces of dimension smaller than $s$. If we use random generators such gaussian or uniform, for true sparse signal $x$, the probability of encountering a ``bad'' $x$ is actually zero.
\end{remark}

% - ``if and only if''?

\subsection{Situation $\rho=1$}
This situation is a little bit tricky: we can not simply assert that every $(P_x)$ has multiple solutions since $A_T$ is square and nonsingular according to the basic \textbf{Assumption \ref{basic_assum}}; the technique in situation $1/2<\rho<1$ does not apply here neither since it is only an inequality. But we can still draw a definitive conclusion with a more careful investigation, shown in below.
\begin{theorem}\label{thm_m=s}
Given a fixed sensing matrix $A\in \setr^{m\times N}$ and we suppose the basic \textbf{Assumption \ref{basic_assum}} holds on $A$. If $\rho=1$, or $s=m$, then every sparse recovery problem $(P_x)$ with $\|x\|_0\leq s$ has at lease two different solutions.
\end{theorem}
\begin{myproof}
We first consider $\|x\|_0=s$. For $x$ with support $S$, we construct another solution $x'$ to $(P_x)$ of the following special form:
\begin{equation}
\supp(x')=S'=(S\backslash \{j\})\cup\{j'\},
\end{equation}
with some $j\in S$ and $j'\in \{1,\ldots,N\}\backslash S$. 

As usual, we let $T=S\cup S'=S\sqcup \{j'\}$ and $|T|=s+1$. We fix $j'\in \{1,\ldots,N\}\backslash S$ and wish to find a $j\in S$ that satisfies our demand. Note that $T$ will not change any more once $j'$ is fixed.

$A_T\in \setr^{m\times (m+1)}$ has a one-dimensional null space and we denote $v$ as a representative:
\begin{equation}
v\neq 0, A_T v=0.
\end{equation}
We now face two situations: there is at least one nonzero entry in $v$ with index lying in $S$ or the only nonzero entry is $v_{j'}$. We claim that the latter situation never happens. If that happens, we have $v_i=0$ for all $i\in S$ and $v_j\neq 0$. This can only lead to $j'$-th column of $A$ equalling to zero, which contradicts with the basic assumption on $A$.

Now we have $j\in S$ and $v_j\neq 0$. Without loss of generality, we assume $v_j=x_j\neq 0$, otherwise we can just scale $v$ to be so. Let 
\begin{equation*}
x'=x-v,
\end{equation*}
then $x'_j = 0$, which means $\|x'\|_0\leq s$ and $x'\neq x$.

If $\|x\|_0<s$, we can artificially complement $\supp(x)$ to have length $s$, and the rest of the proof is exactly the same as above. \qed 
\end{myproof}
\subsection{Summary on uniqueness with respect to oversampling ratio}
If we fix number $s$ of non-zero entries in true sparse signal $x$, the situation gets worse as the number $m$ of measurements decreases.
First, the strongest condition $m\geq 2s$ yields the ``strongest'' uniqueness: the solution to $(P_x)$ is unique for ALL $\|x\|_0\leq s$. Next, as the condition gets weaker ($1/2< s/m<1$), the set of ``bad'' $x$'s becomes a union of some low-dimensional subspaces and as number of measurements $m$ decreases, their dimensions get larger (Theorem \ref{thm_m<2s}). Then $m$ decreases to be the same as $s$ ($s/m=1$) and every $(P_x)$ has multiple solutions (Theorem \ref{thm_m=s}), though we still need to construct a different solution using a different support. Finally, when $s/m>1$, even if we know the support $S$ of true sparse signal $x$, there are still infinitely many solutions to ($P_x$) since sub-matrix $A_S$ itself is a ``fat'' matrix and thus column rank deficient.
%\subsection{Remark on what the constraints ``$A_T(x_T'-x_T'')=0$'' yields, where $T=S'\cup S''$, $|S'|,|S''| \leq s$, $s\leq m<2s$}

\section{Deflation-and-projection-HTP (DP-HTP) algorithm}

In \textbf{Figure \ref{fig:phase_HTPvsDPHTP}}, we observe that the success region of DP-HTP algorithm is significantly larger than that of pure HTP algorithm. In the last section, theoretical analysis shows that it is possible for an algorithm to success even when $1/2<\rho<1$ (the sparse recovery problem ($P_x$) has unique solution for almost all $x$).

In this section, we state deflation-and-projection-HTP algorithm and the rationale behind. Then we make attempts to justify the transition map behaviour of DP-HTP algorithm via ``fast descendent'' condition on true sparse signal.

The HTP algorithm~\cite{foucart2011hard} serves as a critical ingredient in modified algorithm DP-HTP and thus we state it here in \textbf{Algorithm \ref{alg_HTP}} for later use.

\begin{algorithm}[!h]
\caption{Pure HTP algorithm}
\label{alg_HTP}
\KwIn{$A\in \setr^{m\times n}$ with normalized columns, $y\in \setr^n$, sparsity $s$, initial guess $x^{(0)}$ }
\KwOut{$\hat{x}=\mbox{HTP}(A,b,s,x^{(0)})$}
\For{$n=0,1,2,\dots$ until convergence}
{
  $u\rsup{n+1} = x\rsup{n} - \nabla f(x\rsup{n}) = x^{(k)} + A^T(y-Ax\rsup{n})$\;
  $v\rsup{n+1} = \arg\min_{v\in \setr^N}\left\{ \|v-u\rsup{n}\|_2 \big | \|v\|_0 \leq s  \right\} $\;
  $S\rsup{n+1} = \supp(v\rsup{n+1})$\; 
  $x\rsup{n+1} = \arg\min_{z\in \setr^N} \left\{\|y-Az\|_2\big | \supp(z)\subseteq S\rsup{n+1}\right\}$;
  ($x\rsup{n+1}_{S\rsup{n+1}}=A_{S\rsup{n+1}}^\dagger y$, $x\rsup{n+1}_{\overline{S\rsup{n+1}}}=0$)\;
}
$\hat{x}=x\rsup{n+1}$;
\end{algorithm}

\subsection{Statement of DP-HTP algorithm}
If we know part of the solution, we may be able to modify the original problem into an easier one. Deflation is a term for this philosophy. In order to apply this idea, we need to investigate the sparse recovery problem from a different angle: instead of solving $x$ in problem ($P_x$), we only need to find its support $S=\supp(x)$. Now if we know in advance some indices in support $S$, how do we find the rest ones? 

Formally, let $A=[A_1,A_2]$, $x=[x_1^T, x_2^T]^T$, $Ax=A_1x_1+A_2x_2=b$, $x_1\in \setr^{s_1},s_1<s$, and suppose $A_1$ is known. (Note that we do not suppose $x_1$ is known.) Also we assume that all the elements in $x_1$ are non-zero (thus only a few entries in $x_2$ are non-zero since $x$ itself is sparse).

We can not just solve 
\begin{equation}
(P'_{(A_2,x_2)})\quad
  \begin{cases}
  \quad \min_{z\in \setr^{N-s_1}} & \frac{1}{2}\|b-A_2z\|_2^2 \nonumber\\
  \quad \mbox{s.t. } & \|z\|_0\leq s-s_1, \nonumber
  \end{cases}
\end{equation}
since $y\neq A_2 x_2$. Nor can we solve
\begin{equation}
(P''_{(A_2,x_2)})\quad
  \begin{cases}
  \quad \min_{z\in \setr^{N-s_1}} & \frac{1}{2}\|b - A_1\hat{x}_1 - A_2z\|_2^2 \nonumber\\
  \quad \mbox{s.t. } & \|z\|_0\leq s-s_1, \nonumber
  \end{cases}
\end{equation}
where $\hat{x}_1 = A_1^\dagger b$, since $y\neq A_1\hat{x}_1 + A_2x_2$. 

But we can do the following: let $P_1\in \setr^{m\times m}$ be the orthogonal projection onto the orthogonal complement of $\mbox{span}(A_1)$, i.e.,
\begin{equation}
P_1 = I-A_1(A_1^TA_1)^{-1}A_1^T,
\end{equation}
and then multiply $P_1$ on both sides of $A_1x_1+A_2x_2=b$ to get 
\begin{equation}
P_1A_1x_1+P_1A_2x_2 = P_1 b.
\end{equation}
Since $P_1A_1\equiv 0$, we obtain the following size-reduced sparse recovery problem:
\begin{equation}
(P_{(A_2,x_2)})\quad
  \begin{cases}
  \quad \min_{z\in \setr^{N-s_1}} & \frac{1}{2}\|P_1 b -  P_1 A_2z\|_2^2 =  \frac{1}{2}\|b - A_1\hat{x}_1 -P_1 A_2z \|_2^2 \nonumber\\
  \quad \mbox{s.t. } & \|z\|_0\leq s-s_1. \nonumber
  \end{cases}
\end{equation}
Note that $P_1b$ is just the residual after least squares solving on $A_1$. 

The difference between ($P''_{(A_2,x_2)}$) and ($P_{(A_2,x_2)}$) is that the sensing matrix in the latter problem is $P_1A_2$, a projected one, rather than $A_2$. 

Motivated by the above analysis, we design a three-step iteration below:
\begin{itemize}
\item[1] Try HTP on the sparse recovery problem;
\item[2] Select an index based on the result given by HTP;
\item[3] Use projection to shrink the problem into a smaller one.
\end{itemize}
Formally, we write the above iteration steps as the DP-HTP algorithm:
\begin{algorithm}[!ht]
\caption{Deflate-and-Projection HTP algorithm}
\label{alg_DPHTP}
\KwIn{$A\in \setr^{m\times N}$ with normalized columns, $y\in \setr^m$, sparsity $s$}
\KwOut{$\hat{x}=\mbox{DPHTP}(A,b,s)$}
$A^{(0)}=A, b^{(0)}=b, s^{(0)}=s, S^{(0)}=\emptyset$(for index storage) \;
\For{$n=0,1,2,\dots,s-1$}
{
  $x\rsup{n}=\mbox{HTP}(A\rsup{n}, b\rsup{n}, s\rsup{n})$\label{DP_HTP_step_1}\;
  $i\rsup{n}=\arg\max_{i} \{|(x\rsup{n})_i|\}$ \;
  Find $j\rsup{n}\in \{1,\ldots, N\}$ corresponding to $i\rsup{n}$\ ($x\rsup{n} \in \setr^{N-s}$)\;
  $S\rsup{n+1}=S\rsup{n}\cup \{j\rsup{n}\}$ \;
  Write $A\rsup{n}=[A\rsup{n}_1,\ldots, A\rsup{n}_{i\rsup{n}-1},A\rsup{n}_{i\rsup{n}},A\rsup{n}_{i\rsup{n}+1},\ldots, A\rsup{n}_{N-n}]$ and \mbox{\quad\quad\quad} 
  let $\tilde{A}\rsup{n+1}=[A\rsup{n}_1,\ldots, A\rsup{n}_{i\rsup{n}-1},\quad\quad\quad  A\rsup{n}_{i\rsup{n}+1},\ldots, A\rsup{n}_{N-n}]$\;
  $A\rsup{n+1}=\left[ I- A\rsup{n}_{i\rsup{n}}(A^{(n)T}_{i\rsup{n}}A\rsup{n}_{i\rsup{n}})^{-1}A^{(n)T}_{i\rsup{n}}\right]\tilde{A}\rsup{n+1}$; (project columns in $\tilde{A}\rsup{n+1}$ onto the orthogonal complement of $\mbox{span}\{A\rsup{n}_{i\rsup{n}}\}$)\
  $b\rsup{n+1} = \left[ I- A\rsup{n}_{i\rsup{n}}(A^{(n)T}_{i\rsup{n}}A\rsup{n}_{i\rsup{n}})^{-1}A^{(n)T}_{i\rsup{n}}\right]b\rsup{n} \in \setr^m$ \;
  $s\rsup{n+1}=s\rsup{n}-1$\;
}
$\hat{x}_{S^{(s)}}=A^\dagger_{S^{(s)}} y$, $\hat{x}_{\overline{S^{(s)}}}$=0;
\end{algorithm}

\subsection{Justification of DP-HTP algorithm}
For terminology simplification, we denote from now on true sparse signal by $x^\ast, \|x^\ast\|\leq s$, its support $S^\ast=\supp(x^\ast)$, and measurement vector $b = Ax^\ast$. The DP-HTP algorithm calls HTP algorithm in each iteration, picks up the largest entry in the result and then projects the problem into an ``one-order-smaller''(using smaller is better?) one. 

There are two conditions for DP-HTP algorithm to work: that the HTP step always returns a result (though it may not be the correct solution) and that the index-selection step always picks up a correct index. Although our numerical experiments show that the iterate in HTP algorithm never become periodic in practice, there seem to be no theoretical characterisation for the situation here. Thus we leave the former condition as an assumption. For the latter one, out attempt is to impose assumptions on both the sensing matrix and the distribution of the non-zero entries in the true sparse signal.

For the sensing matrix part, we have the following theoretical result:

\begin{theorem}
The $(k-1)$-th order restricted isometry constant of $A\rsup{n+1}$ can be controlled by $k$-th order RIC of $A\rsup{n}$, i.e., if there exists $\delta_k\rsup{n}<1$ s.t., 
\begin{equation}
(1-\delta_k\rsup{n})x^Tx\leq x^TA^{(n)T}A\rsup{n}x\leq (1+\delta_k\rsup{n})x^Tx, \quad \forall \|x\|_0\leq k,
\end{equation}
then there exists $\delta_{k-1}\rsup{n+1}<1$ such that
\begin{equation}
(1-\delta_{k-1}\rsup{n+1})\tilde{x}^T\tilde{x}\leq \tilde{x}^TA^{(n+1)T}A\rsup{n+1}\tilde{x}\leq (1+\delta_{k-1}\rsup{n+1})\tilde{x}^T\tilde{x},\quad \forall \|\tilde{x}\|_0\leq k-1,
\end{equation}
and 
\begin{equation}
\delta_{k-1}\rsup{n+1} \leq \delta_{k}\rsup{n}.
\end{equation}
\end{theorem}
\begin{myproof}
Verifying restricted isometry property of $A\rsup{n+1}$ involves estimating the singular values of sub-matrices with $k-1$ columns.

Choose any $k-1$ columns of $A\rsup{n+1}$, denoted by $\tilde{B}_2\in \setr^{m\times (k-1)}$. Find the corresponding columns in $A\rsup{n}$, and name them as $B_2\in \setr^{m\times (k-1)}$. Through the definition of $A\rsup{n+1}$, we have $\tilde{B}_2 = P\rsup{n}B_2$. Let $B_1 = A_{i\rsup{n}}\rsup{n}$ and $B=[B_1,B_2]$. We are to bound the singular values of $\tilde{B_2}$ by those of $B$.

Let $[U_1, U_2]$ be a set of orthogonal basis of $\mbox{span}(B)$ with $U_1=B_1/\|B_1\|_2$. Since $P\rsup{n}$ is the projection onto the orthogonal complement of $\mbox{span}(B_1)$, 
\begin{equation}
\tilde{B}_2 = P\rsup{n}B_2=(I-U_1U_1^T)B_2=U_2U_2^TB_2.
\end{equation}
Write
\begin{equation}
B=UU^TB=
\left[
\begin{matrix}
U_1 & U_2
\end{matrix}
\right]
\left[
\begin{matrix}
 U_1^T B_1 & U_1^TB_2 \\ 
 0         & U_2^TB_2
 \end{matrix}
 \right],
\end{equation}
and we have 
\begin{itemize}
\item The $k$ singular values of $B$ are exactly the same as those of 
  \begin{equation}
    \left[
		\begin{matrix}
		 U_1^T B_1 & U_1^TB_2 \\ 
		 0         & U_2^TB_2
		\end{matrix}
		\right],
  \end{equation}
\item The $(k-1)$ singular values of $\tilde{B}_2 = U_2U_2^T B_2$ are exactly the same as those of  $U_2^TB_2$.
\end{itemize}
By Cauchy interlacing theorem~\cite{horn2012matrix}, the ordered singular values of $U_2^TB_2$ can be inserted into the order singular values of $B$, which implies
\begin{equation}
\sqrt{1-\delta_k\rsup{n}}\leq \sigma(B)\leq \sqrt{1+\delta_k\rsup{n}} \Rightarrow \sqrt{1-\delta_k\rsup{n}}\leq \sigma(U_2^TB_2)\leq \sqrt{1+\delta_k\rsup{n}}
\end{equation}
and thus implies $\delta_{k-1}\rsup{n+1} \leq \delta_{k-1}\rsup{n+1}$. \qed
\end{myproof}

\begin{theorem} \label{RIC_of_A}
The $(k-n)$-th order RIC $\delta_{k-n}\rsup{n}$ of $A\rsup{n}$ can be controlled by the $k$-order RIC $\delta_k$ of $A$, i.e., $\delta_{k-n}\rsup{n} \leq \delta_k$.
\end{theorem}

\textbf{Theorem \ref{RIC_of_A}} asserts that the smaller
sensing matrix is never worse than the previous one in the sense of restricted isometry property. Our next mission is to derive conditions that ensure the index of the largest entry in the result given by HTP algorithm corresponds to some index lying in the support of true solution $\supp(x^\ast)$. Note that it is impossible to develop any algorithm that can work for ALL $x^\ast, \|x^\ast\|\leq s$ when $m<2s$, since the sparse recovery problem itself may have multiple solutions then (see \textbf{Theorem \ref{thm_m<2s}}). Thus it is natural to assume constraints on the true sparse signal $x^\ast$. Our approach is to analysis one(the first) iteration in DP-HTP algorithm and see whether it is compatible with mathematical induction.

Let $x=\mbox{HTP}(A,b,s)$ be the result given by HTP algorithm (the first iteration in DP-HTP algorithm), and $S=\supp(x)$ be its support. Furthermore, we give three conditions below:

\begin{assumption}\label{cond_1}
  $A$ has $(s+2)$-th order restricted isometry constant $\delta=\delta_{s+2}<1$. Its detailed bound will be derived later.
\end{assumption} 
\begin{assumption}\label{cond_2}
  HTP algorithm on $(A,b,s)$ converges (may not converge to true sparse signal $x^\ast$).
\end{assumption}
\begin{assumption}\label{cond_3}
  The largest entry in $x^\ast$ is dominant by the following definition:
\begin{equation}\label{eq_dominant}
  l = \arg\max_i\{|x^\ast_i|\| i=1,\leq,N\}, \mbox{and } |x^\ast_l|\geq \gamma \|x^\ast_{S^\ast\backslash\{l\}}\|_2.
\end{equation}
\end{assumption}
Our aim is to derive some condition involving restricted isometry constant $\delta$ and dominance factor $\gamma$ that can ensure $l\in S^\ast$. 

Let $u = x+A^T(b-Ax)$. According to \textbf{Assumption \ref{cond_2}}, we only need to ensure
\begin{equation}
    |u|_l > |u|_j,\ \forall j\in S\backslash S^\ast.
\end{equation}
By the definition of HTP algorithm, the nonzero entries $x_S$ are solved via least squares:
\begin{equation}
x_S=(A^T_S A_S)^{-1}A_S^Tb = (A^T_S A_S)^{-1}A_S^TA_{S^\ast}x^\ast_{S^\ast}.
\end{equation}
Also we denote the orthogonal project onto the orthogonal complement of $\mathrm{span}(A_S)$ by
\begin{equation}
P_S \equiv I-A_S(A^T_S A_S)^{-1}A_S^T
\end{equation}
for latter use.
\begin{proposition}\label{prop_prod_1}
    Let $B=[B_1, B_2]$, $B_1\in \setr^{m\times s}$, $B_2\in \setr^{m\times 1}$ ($m\geq s+1$), and orthogonal projection $P=I-B_1(B_1^TB_1)^{-1}B_1^T$. Suppose $\|B^TB-I\|_2\leq \delta$, or equivalently, singular values of $B$ are bounded by $[\sqrt{1-\delta}, \sqrt{1+\delta}]$. Then $1-\delta \leq \langle B_2, PB_2\rangle \leq 1+\delta$.
\end{proposition}
\begin{myproof}
    Using Gram-Schmidt orthogonalization, we can assume that the columns of $U\in \setr^{m\times (s+1)}$ are the orthogonal basis of $\mathrm{span}(B)$, $U=[U_1, U_2]$, and columns of $U_1\in \setr^{m\times s}$ are the orthogonal basis of $B_1$, which yields $U_2^T B_1=0$. Write
    \begin{equation}
        B = UU^TB = 
        \left[
        \begin{matrix}
            U_1& U_2 
        \end{matrix}
        \right]
        \left[
        \begin{matrix}
            U_1^TB_1& U^T_1 B_2  \\ 
            0       & U^T_2 B_2 
        \end{matrix} 
        \right].
    \end{equation}
    $U^T_2 B_2$ (a scalar) is the $2$-by-$2$ block, and by the Cauchy interlacing theorem, it can be bounded by the singular values of $B$: $\sqrt{1-\delta}\leq U^T_2 B_2 \leq  \sqrt{1+\delta}$. Thus $\langle B_2, PB_2 \rangle = B_2^T U_2 U_2^T B_2 \in [1-\delta, 1+\delta]$. \qed
\end{myproof}

\begin{proposition}\label{prop_prod_2}
    Let $B=[B_1, B_2, B_3]$, $B_1\in \setr^{m\times s}$, $B_2, B_3\in \setr^{m\times 1}$ ($m\geq s+2$), and orthogonal projection $P=I-B_1(B_1^TB_1)^{-1}B_1^T$. Suppose $\|B^TB-I\|_2\leq \delta$, or equivalently, singular values of $B$ are bounded by $[\sqrt{1-\delta}, \sqrt{1+\delta}]$. Then $|\langle B_3, PB_2\rangle|\leq \delta$.
\end{proposition}
\begin{myproof}
    Similar to the proof in proposition \ref{prop_prod_1}, we assume there exists orthonormal $U=[U_1, U_{23}]\in \setr^{m\times(s+1)}$, $U_1\in \setr^{m\times s}$, $\mathrm{span}(U)=\mathrm{span}(B)$,  and $\mathrm{span}(U_1)=\mathrm{span}(B_1)$. Write identity 
    \begin{equation}
        B = UU^TB = [U_1, U_{23}]
        \left[
        \begin{matrix}
            U_1^T B_1& U_1^TB_{23}  \\ 
            0        & U_{23}^TB_{23} 
        \end{matrix} 
        \right].
    \end{equation}
    $C=U_{23}^TB_{23} \in \setr^{2\times 2}$ is the 2-by-2 block and thus has bound $\|C^TC-I\|_2\leq \delta$. Again we separate blocks of $C^TC-I$ by
    \begin{equation}
        C^TC-I = 
        \left[
        \begin{matrix}
            B_2^TU_{23}U_{23}^TB_2-1& B_2^TU_{23}U_{23}^TB_3  \\ 
            B_3^TU_{23}U_{23}^TB_2  & B_3^TU_{23}U_{23}^TB_3-1
        \end{matrix} 
        \right],
    \end{equation}
    then 
    \begin{equation}
        |\langle B_3, PB_2\rangle| = \|B_3^TU_{23}U_{23}^TB_2\|_2\leq \|C^TC-I\|_2\leq \delta,
    \end{equation}
    since the 2-norm of a sub-matrix is bounded by that of the full matrix.\qed
\end{myproof}

Suppose $l\in S^\ast\backslash S$. We are to derive a constraint involving $\delta$ and $\gamma$ that will yield $|u|_l>|u|_j,\ \forall j\in \overline{S\cup S^\ast}$, which will in turn lead to contradiction. To separate the term containing $x^\ast_l$, we write
\begin{align}
u_l &=\langle A_l, b-Ax \rangle,\  (\mbox{since } x_l=0)\\
    &= \langle A_l, A_{S^\ast}x^\ast_{S^\ast} -A_S(A^T_S A_S)^{-1}A_S^TA_{S^\ast}x^\ast_{S^\ast}\rangle\\
    &=\langle A_l, P_S(A_{S^\ast\cap S} x^\ast_{S^\ast\cap S}+A_{S^\ast\backslash S}x^\ast_{S^\ast\backslash S})\rangle\\
    &=\langle A_l, P_S A_l\rangle x^\ast_l + \langle A_l, P_S A_{S^\ast \backslash (S\cup \{l\})}x^\ast_{S^\ast \backslash (S\cup \{l\})}\rangle.\label{ul_terms}
\end{align}
By \textbf{Proposition \ref{prop_prod_1}}, term $\langle A_l, P_S A_l\rangle$ in (\ref{ul_terms}) can bounded from below by
\begin{equation}
|\langle A_l, P_S A_l\rangle|\geq 1-\delta,
\end{equation}
then $|u_l|$ has a lower bound
\begin{align}
|u_l|&\geq (1-\delta)|x^\ast_l| - \|A_l\|_2\|P_S\|_2\|A_{S^\ast \backslash (S\cup \{l\})}\|_2\|x^\ast_{S^\ast \backslash (S\cup \{l\})}\|_2\\
     &\geq (1-\delta)|x^\ast_l|-\sqrt{1+\delta_1(A)}\sqrt{1+\delta_{s-1}(A)}\|x^\ast_{S^\ast \backslash (S\cup \{l\})}\|_2\\
     &\geq
     \left(
     1-\delta-\frac{1+\delta}{\gamma}\frac{\|x^\ast_{S^\ast \backslash (S\cup \{l\})}\|_2}{\|x^\ast_{S^\ast \backslash \{l\}}\|_2}
     \right)
     |x^\ast_l|.\label{ineq_u_1}
\end{align}
The last inequality is derived via \textbf{Assumption \ref{cond_3}}.

Using similar techniques and proposition \ref{prop_prod_2}, we derive the upper bound of $|u_j|$, $j\in \overline{S\cup S^\ast}$:
\begin{align}
  |u_j| &= |\langle A_j, P_SA_lx^\ast_l\rangle+\langle A_j, P_SA_{S^\ast\backslash (S\cup \{l\})}x^\ast_{S^\ast\backslash (S\cup \{l\})} \rangle|\\
        &\leq \left(\delta + \frac{1+\delta}{\gamma}\frac{\|x^\ast_{S^\ast \backslash (S\cup \{l\})}\|_2}{\|x^\ast_{S^\ast \backslash \{l\}}\|_2}\right)|x^\ast_l|.\label{ineq_u_2}
\end{align}
If the dominance factor $\gamma$ satisfies
\begin{equation}\label{ineq_gamma_1}
  \gamma > \frac{2(1+\delta)}{1-2\delta}\frac{\|x^\ast_{S^\ast \backslash (S\cup \{l\})}\|_2}{\|x^\ast_{S^\ast \backslash \{l\}}\|_2},
\end{equation}
then $|u_l|>|u_j|$, $\forall j\in \overline{S\cup S^\ast}$, implying contradiction on the supposition $l\in S^\ast\backslash S$.

Now we turn to the case when $l\in S\cap S^\ast$, and we are to find conditions to assure $|u_l|>|u_i|$, $\forall i\in S\backslash S^\ast$. Similarly, we derive lower bound for $u_l$ and upper bound for $u_i$. On one hand, 
\begin{align}
|u_l|=|x_l|&=|\langle e_l, (A_S^TA_S)^{-1}A_S^T(A_{S^\ast\cap S}x^\ast_{S^\ast\cap S}+ A_{S^\ast\backslash S}x^\ast_{S^\ast\backslash S})\rangle|\\
     &=|x_l^\ast+\langle e_l,A_{S^\ast\backslash S}x^\ast_{S^\ast\backslash S}\rangle|\\
     &\geq |x_l^\ast| -\|(A_S^TA_S)^{-1}A_S^T\|_2 \|A_{S^\ast\backslash S}\|_2\|x^\ast_{S^\ast\backslash S}\|_2.
\end{align}
Considering the fact that the singular values of $(A_S^TA_S)^{-1}A_S^T$ are the reciprocal of those of $A_S$, we have 
\begin{equation}\label{ineq_u_3}
|u_l|\geq
\left(
1-\frac{1}{\gamma}
\sqrt{\frac{1+\delta}{1-\delta}}
\frac{\|x^\ast_{S^\ast\backslash S}\|_2}{\|x^\ast_{S^\ast\backslash \{l\}}\|_2}
\right)|x_l^\ast|.
\end{equation}
On the other hand, $\forall i\in S\backslash S^\ast$, we have similarly,
\begin{align}
  |u_i|=|x_i|&=|\langle e_i, (A_S^TA_S)^{-1}A_S^TA_{S^\ast\backslash S}x^\ast_{S^\ast\backslash S}\rangle|\\
  &\geq \frac{1}{\gamma}\sqrt{\frac{1+\delta}{1-\delta}}
  \frac{\|x^\ast_{S^\ast\backslash S}\|_2}{\|x^\ast_{S^\ast\backslash \{l\}}\|_2}
  |x_l^\ast|.\label{ineq_u_4}
\end{align}
Again, as long as the dominance factor $\gamma$ satisfies
\begin{equation}\label{ineq_gamma_2}
  \gamma>2\sqrt{\frac{1+\delta}{1-\delta}}
  \frac{\|x^\ast_{S^\ast\backslash S}\|_2}{\|x^\ast_{S^\ast\backslash \{l\}}\|_2},
\end{equation}
then $|u_l|>|u_i|$, $\forall i\in S\backslash S^\ast$. Combining inequalities (\ref{ineq_gamma_1}) and (\ref{ineq_gamma_2}) together, we conclude that the index of the largest entry in $x$, though may not be $l$, is assured to lie in the support $S^\ast$ of the true sparse signal, if the following inequality holds:
\begin{align}
  \gamma&>\max\left\{
  \frac{2(1+\delta)}{1-2\delta}\frac{\|x^\ast_{S^\ast \backslash (S\cup \{l\})}\|_2}{\|x^\ast_{S^\ast \backslash \{l\}}\|_2},
  2\sqrt{\frac{1+\delta}{1-\delta}}
  \frac{\|x^\ast_{S^\ast\backslash S}\|_2}{\|x^\ast_{S^\ast\backslash \{l\}}\|_2}
  \right\} \nonumber\\
        &\geq \frac{2(1+\delta)}{1-2\delta}\frac{\|x^\ast_{S^\ast \backslash (S\cup \{l\})}\|_2}{\|x^\ast_{S^\ast \backslash \{l\}}\|_2}.\label{ineq_gamma}
\end{align}

\begin{remark}
  If we simply relax the factor $\|x^\ast_{S^\ast \backslash (S\cup \{l\})}\|_2/\|x^\ast_{S^\ast \backslash \{l\}}\|_2$ to be one, the dominance factor $\gamma$ will need to be larger than $2(1+\delta)/(1-2\delta)\geq 2$ to guarantee validity, according to (\ref{ineq_gamma}). This makes the assumption \ref{cond_3} too strong to satisfy. However, if the support $S$ given by HTP algorithm already contains a large part of the true support $S^\ast$, the factor $\|x^\ast_{S^\ast \backslash (S\cup \{l\})}\|_2/\|x^\ast_{S^\ast \backslash \{l\}}\|_2$ can be expected to be far smaller than one, making it possible for a random sparse signal with Gaussian distribution to satisfy \textbf{Assumption \ref{cond_3}}. The heuristic observations on factor $\|x^\ast_{S^\ast \backslash (S\cup \{l\})}\|_2/\|x^\ast_{S^\ast \backslash \{l\}}\|_2$ are discussed in numerical examples in the next section. (need completion)
\end{remark}

\begin{remark}Comparison with OMP-type strategy. The DP-HTP algorithm utilizes the result given by HTP algorithm. One may wonder whether HTP is necessary: after all, each HTP requires 10-20 iterations and each iteration involves a least squires problem solving. Here, we discuss a simplified version: using an Orthogonal Matching Pursuit (OMP)-type iteration instead of HTP iteration in each DP-HTP step (see step \ref{DP_HTP_step_1} in \textbf{algorithm \ref{alg_DPHTP}}). Since OMP uses $A^TAb$ as the support indicator, we estimate
\begin{align}
  |u_l| &=|\langle A_l, A_l x_l^\ast \rangle + \langle A_l, A_{S^\ast \backslash \{l\}}x^\ast_{S^\ast\backslash \{l\}}\rangle|\nonumber\\
        &\geq (1-\delta)|x_l^\ast| - \delta\|x^\ast_{S^\ast\backslash \{l\}}\|_2\\
        &\geq (1-\delta-\frac{\delta}{\gamma}),
\end{align}
and $\forall j\neq l$,
\begin{align}
  |u_j| &= |\langle A_j, A_l x^\ast_l\rangle + \langle A_j, A_{S^\ast\backslash \{l,j\}}x^\ast_{S^\ast\backslash \{l,j\}}  \rangle + \langle A_j, A_j x^\ast_j \rangle|\\
        &\leq \delta |x^\ast_l| + \delta \|x^\ast_{S^\ast\backslash \{l,j\}}\|_2 + (1+\delta) |x^\ast_j|\\
        &\leq \delta |x^\ast_l| + \sqrt{\delta^2+(1+\delta)^2} \|x^\ast_{S^\ast\backslash \{l\}}\|_2\\
        &\leq \left(\delta + \frac{1}{\gamma}\sqrt{\delta^2+(1+\delta)^2}\right)|x^\ast_l|.
\end{align}
Thus the condition for ensuring $|u_l|>|u_j|$ should be
\begin{equation}\label{ineq_gamma_OMP}
  \gamma > \frac{\delta+\sqrt{\delta^2+(1+\delta)^2}}{1-2\delta}.
\end{equation}
The fundamental difference between (\ref{ineq_gamma}) and (\ref{ineq_gamma_OMP}) is the factor 
 \begin{equation}
   \|x^\ast_{S^\ast \backslash (S\cup \{l\})}\|_2/\|x^\ast_{S^\ast \backslash \{l\}}\|_2,
 \end{equation} 
 which is essential for \textbf{Assumption \ref{cond_3}} to be satisfied on randomly generated true signal $x^\ast$, with, for example, Gaussian distribution.
\end{remark}

\begin{remark}
  The effectiveness of the DP-HTP algorithm heavily relies on the result given by the HTP algorithm in every step, as can be seen from the definition of DP-HTP and the analysis above.
\end{remark}

\section{Numerical Experiments}
\subsection{Gaussian Setting}
\begin{figure}[bh]
	\includegraphics[scale=0.27]{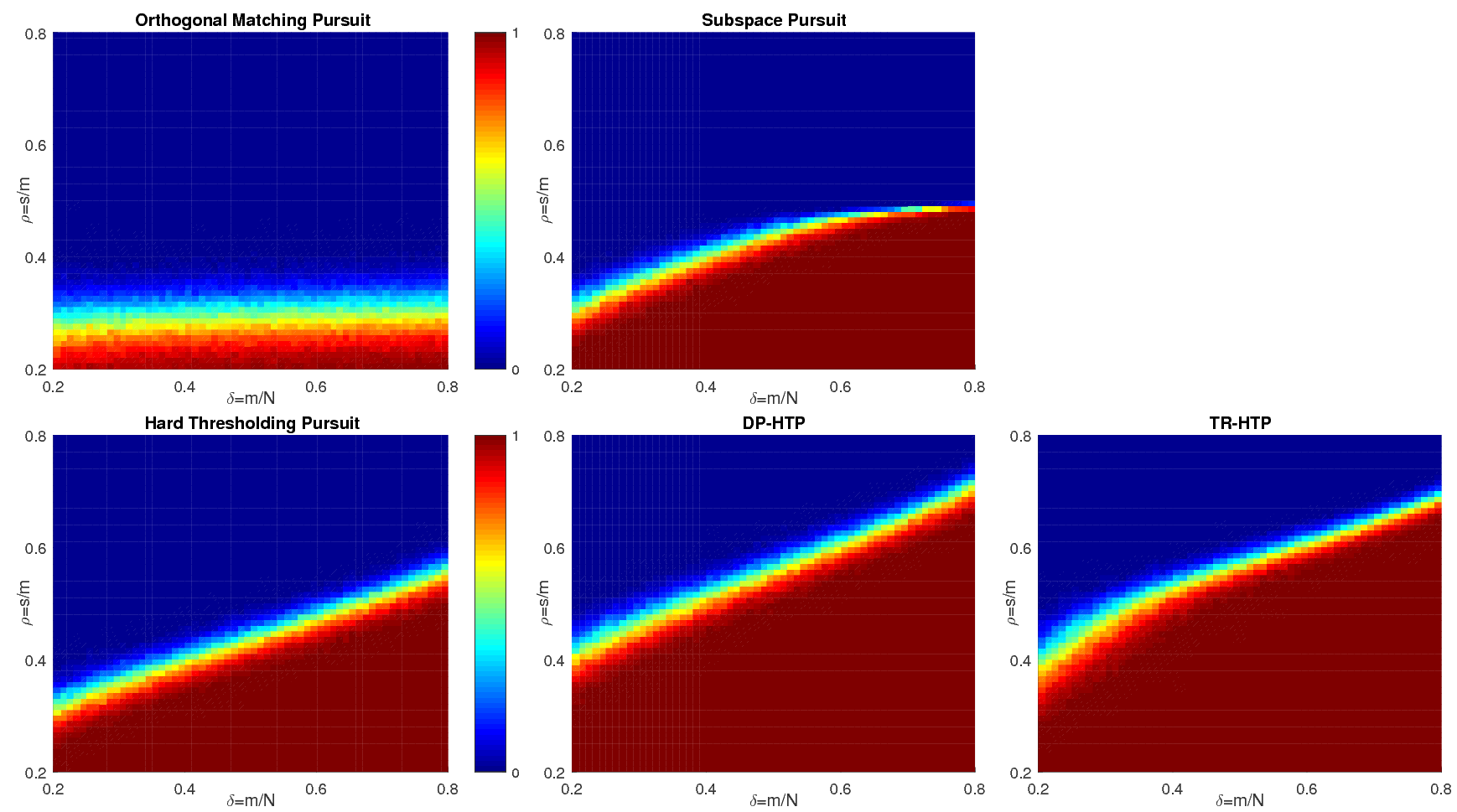}
	\caption{phase transition map for Gaussian setting}
	\label{fig:phase_for_Gaussian}
\end{figure}

\begin{figure}[th]
	\includegraphics[scale=0.36]{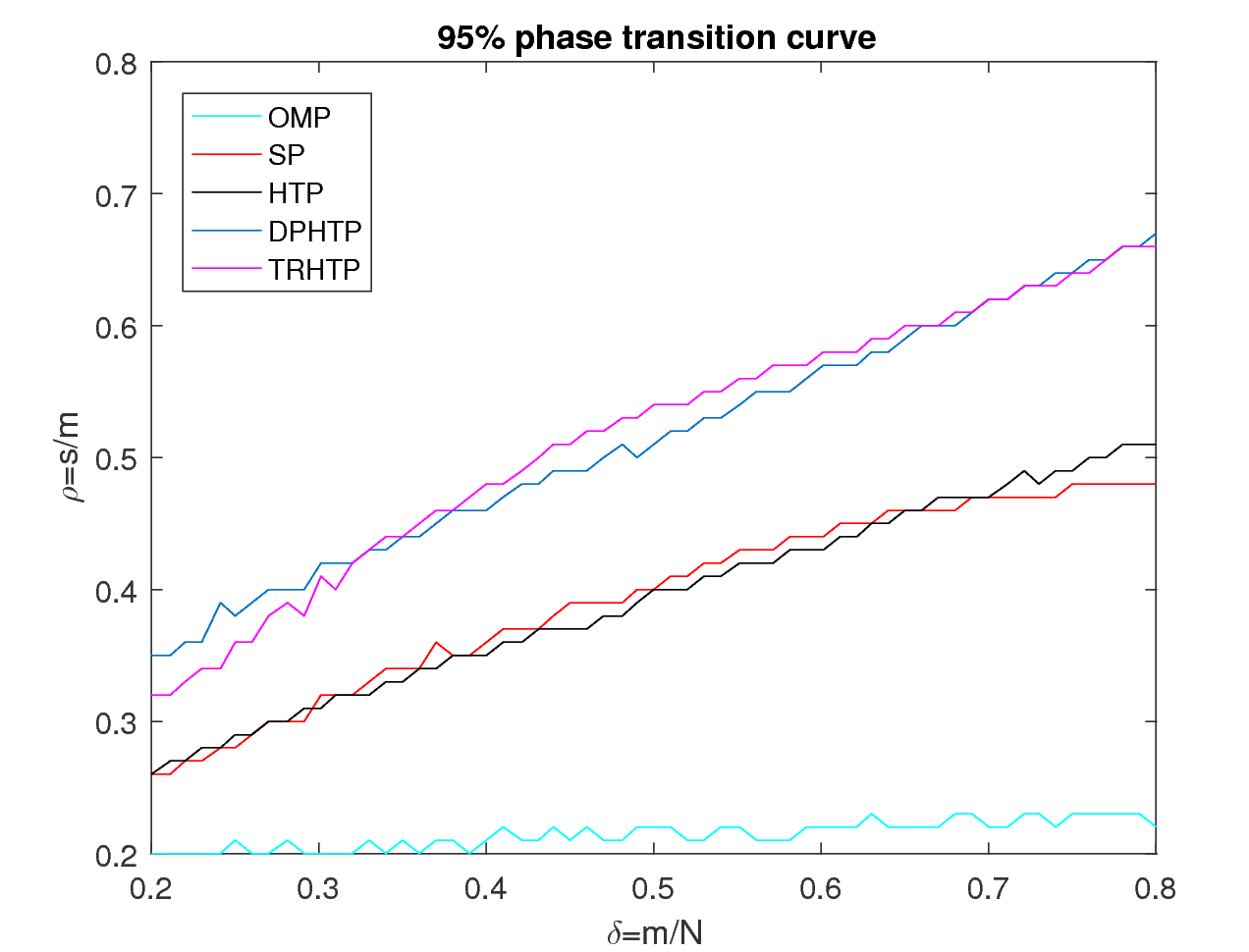}
	\caption{95\% success line for Gaussian setting}
	\label{fig:phase_line_for_Gaussian_95}
\end{figure}

We set the experiment as follows: 
\begin{itemize}
	\item Signal length $N=800$. Number of measurements $m$ and number of nonzeros $s $ are set such that both oversampling ratio $\rho=s/m$ and undersampling ratio $\delta=m/N$ vary from $0.2$ to $0.8$.
	\item Both the entries of sensing matrix $A$ and the nonzero entries of true solution $x$ are $i.i.d$ Gaussians with unit standard variation.
	\item The experient is repeated 1000 times independently.
	\item The algorithms in comparison are OMP, Subspace Pursuit and HTP, as well as modified HTP algorithms proposed in this paper.
\end{itemize}
The overall phase transition map for the five algorithms in comparison can be found in Fig.\ref{fig:phase_for_Gaussian} and Fig.\ref{fig:phase_line_for_Gaussian_95} shows a detailed comparison for 95\% success rate performance. The two modified HTP algorithms TR-HTP and DP-HTP perform far better on phase transition behavior than the rest three greedy algorithms do. Also when undersampling ratio $\delta$ is large (close to $0.8$), TR-HTP and DP-HTP manage to solve $P_x$ with oversampling ratio $\rho$ much larger than $0.5$.

% bibtex
\bibliography{HTPKai.bib}
\bibliographystyle{plain}
\end{document}